\documentclass[journal]{IEEEtran}
\usepackage{graphicx}
\usepackage[normalem]{ulem}
\usepackage{amssymb}
\usepackage{amsmath}
\usepackage{amstext}
\usepackage{amsfonts}
\usepackage[latin1]{inputenc}
\usepackage[normalem]{ulem}
\usepackage[brazil]{babel}
\usepackage{epstopdf}
\usepackage{amsmath}
\usepackage{caption}
\usepackage{hyperref}

\newtheorem{theorem}{Theorem}

\newtheorem{proposition}{Proposition}

\begin{document}
\newcommand{\doi}{doi:}
\renewcommand*{\doi}[1]{\href{http://dx.doi.org/#1}{doi: #1}}
\renewcommand*\abstractname{Abstract}
\renewcommand{\figurename}{Figure}
\renewcommand{\tablename}{Table}
\renewcommand{\refname}{References}

\title{A Kotel'nikov Representation for Wavelets}

\author{{H. M. de Oliveira}, {R. J. Cintra}, {R. C. de Oliveira}\\
Federal University of Pernambuco, UFPE, \\
Statistics Department, CCEN, Recife, Brazil.\\
State University of Amazon, UEA,\\
Computer Engineering Department, Manaus, Brazil.}
\date{2017}

\maketitle

\begin{abstract}

\textbf{This paper presents a wavelet representation using baseband signals, by exploiting Kotel'nikov results. Details of how to obtain the processes of envelope and phase at low frequency are shown. The archetypal interpretation of wavelets as an analysis with a filter bank of constant quality factor is revisited on these bases. It is shown that if the wavelet spectral support is limited into the band $[f_m,f_M]$, then an orthogonal analysis is guaranteed provided that $f_M \leq 3f_m$, a quite simple result, but that invokes some parallel with the Nyquist rate. Nevertheless, in cases of orthogonal wavelets whose spectrum  does not verify this condition, it is shown how to construct an ``equivalent'' filter bank with no spectral overlapping.} 
\end{abstract}

\providecommand{\keywords}{\textbf{\textit{Index terms---}}}

\section{Introduction}

Wavelet and Mallat's multiresolution analysis have become well-established tools in signal analysis, particularly to provide more efficient processing, and as a consequence of the growing availability of techniques \cite{deO}, \cite{Mallat}. One usual interpretation to introduce them is the $Q$-constant filter bank analysis, easy to understand \cite{Vetterli}. Moreover, the bandpass behavior of the wavelet is one of its widely recognized features \cite{Mallat}. On the other side, in telecommunication systems and noise models \cite{Vetterli}, \cite{Haykin_Moher}, perhaps the most relevant representation (and certainly the most used one) is the bandpass representation  \cite{Taub}, involving a modulated baseband signal. As the condition of admissibility of wavelets, $\psi(t) \leftrightarrow \Psi(w)$, imposes a zero in the spectrum origin, i.e., $\Psi(0)=0$, these signals are known to be bandpass. How to apply the usual bandpass representations to wavelets? This is one of the focuses of this paper. In this framework, it rescues the following theorem by Kotel'nikov.\\	
\begin{theorem}
(Kotel'nikov, 1933). \textit{If a signal $f (t)$ has a spectrum confined in the band $[w_1, w_2]$, then there is a representation $f(t)=g(t).cos \left [ \frac{w_1+w_2}{2}t+\theta (t) \right ]$ 
being $g(t)$ and $\theta (t)$ low-frequency processes, limited in $\left [0,\frac{w_2-w_1}{2} \right ].$\\}
\end{theorem}
This result was shown in the pioneering paper demonstrating the sampling theorem (prior to Shannon theorem \cite{Kotelnikov}, \cite{Bissell}). An elegant proof is found in Theorem 4 of \cite{Haykin_Moher}. The paper is organized as follows. Initially, it is introduced in Section II, a bandpass representation to wavelets, showing how to get the envelope and the phase signal. This is applied to a few known wavelets, illustrating such a representation. In Section III, we investigate the effects of scaling (daughter wavelet generation), connecting directly and naturally to the wavelet analysis with filter bank. The imposition of a condition of no spectral overlap of the filters leads to a new condition for orthogonal analysis. Again, it shows that the condition is checked to simple continuous orthogonal wavelet, including Shannon wavelet. For wavelets with asymmetry in spectrum, including Meyer wavelet \cite{Burrus}, \cite{Meyer}, Daubechies (e.g. db4 \cite{Mallat}), wavelet ``de Oliveira'' \cite{deO3}, the analysis is also done. The findings are finally presented in Section IV.
\section{Bandpass Representation to Wavelets}
\begin{figure}
\centering
\includegraphics[scale=0.48]{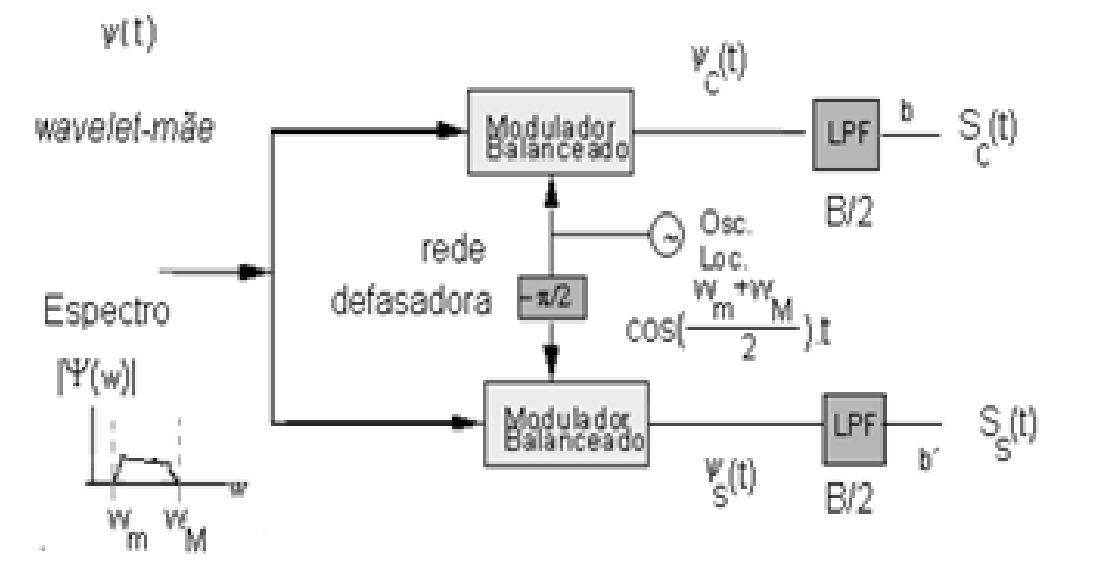}
\captionof{figure}{Decomposition of a wavelet contained in the spectral range $[w_m,w_M]$ in their low-frequency components in phase and in quadrature. The value of $B$ corresponds to the band wavelet, $2\pi B=w_M-w_m.$}
\label{fig:fig1}
\end{figure}

\begin{figure}
\centering
\includegraphics[scale=0.5]{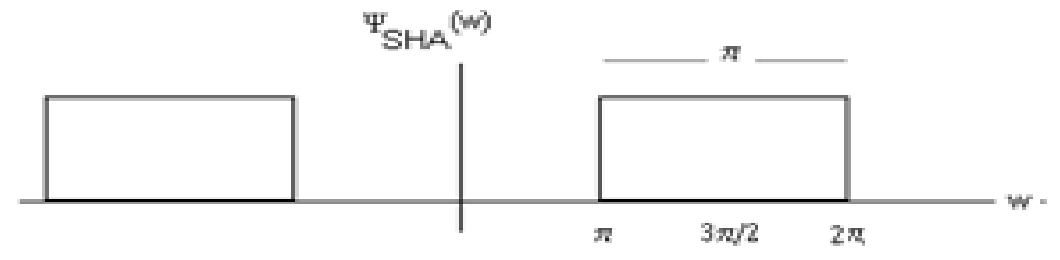}
\captionof{figure}{Spectrum of wavelet Shannon, watching the confinement in the spectral range $[\pi/2,5 \pi/2]$, with center frequency $w_0=3 \pi /2$ and band $2 \pi B=\pi$, i.e. 0.5 Hz.}
\label{fig:fig2}
\end{figure}

\begin{figure}
\centering
\includegraphics[scale=0.55]{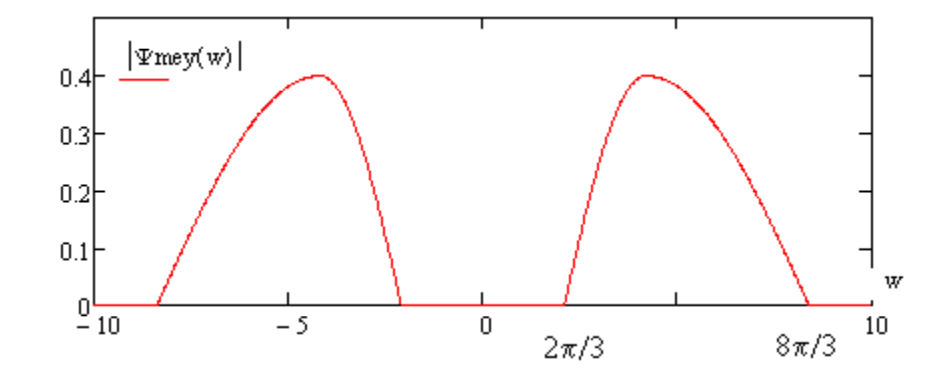}
\captionof{figure}{Meyer wavelet spectrum, observing the confinement in the spectral range $[2\pi/3 8, \pi/3]$, with center frequency $5\pi/2$, but $w_0=6 \pi/2$ and bandwidth $2\pi B=2 \pi$, i.e. 1 Hz.}
\label{fig:fig3}
\end{figure}

\begin{figure}
\centering
\includegraphics[scale=0.45]{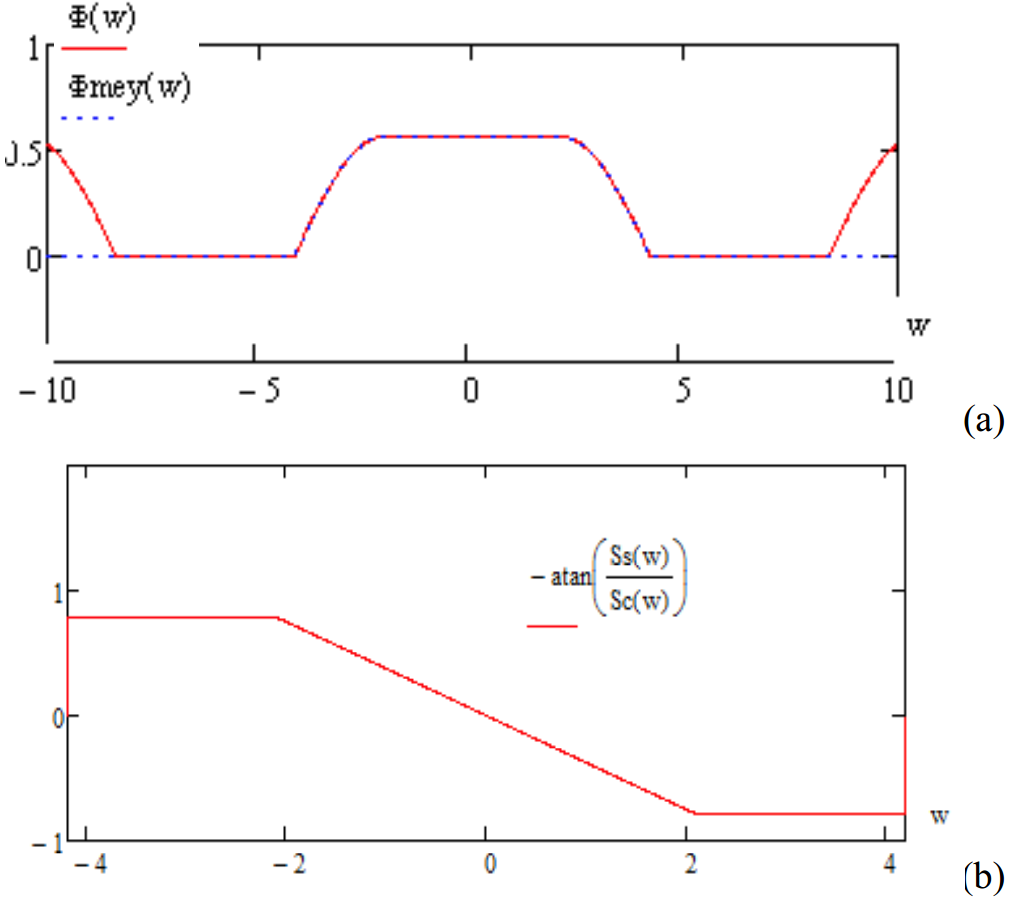}
\captionof{figure}{Spectrum of the baseband component for the Meyer wavelet, which corresponds to the envelope of the wavelet. The spectrum is zero for $|w| \geq \pi$, featuring a limited bandwidth signal at 0.5 Hz b) representation of the phase of the spectrum $\Phi(w)$ setting in terms of $S_c$ and $S_s$ components, cf.(12).}
\label{fig:fig4}
\end{figure}

\begin{figure}
\centering
\includegraphics[scale=0.45]{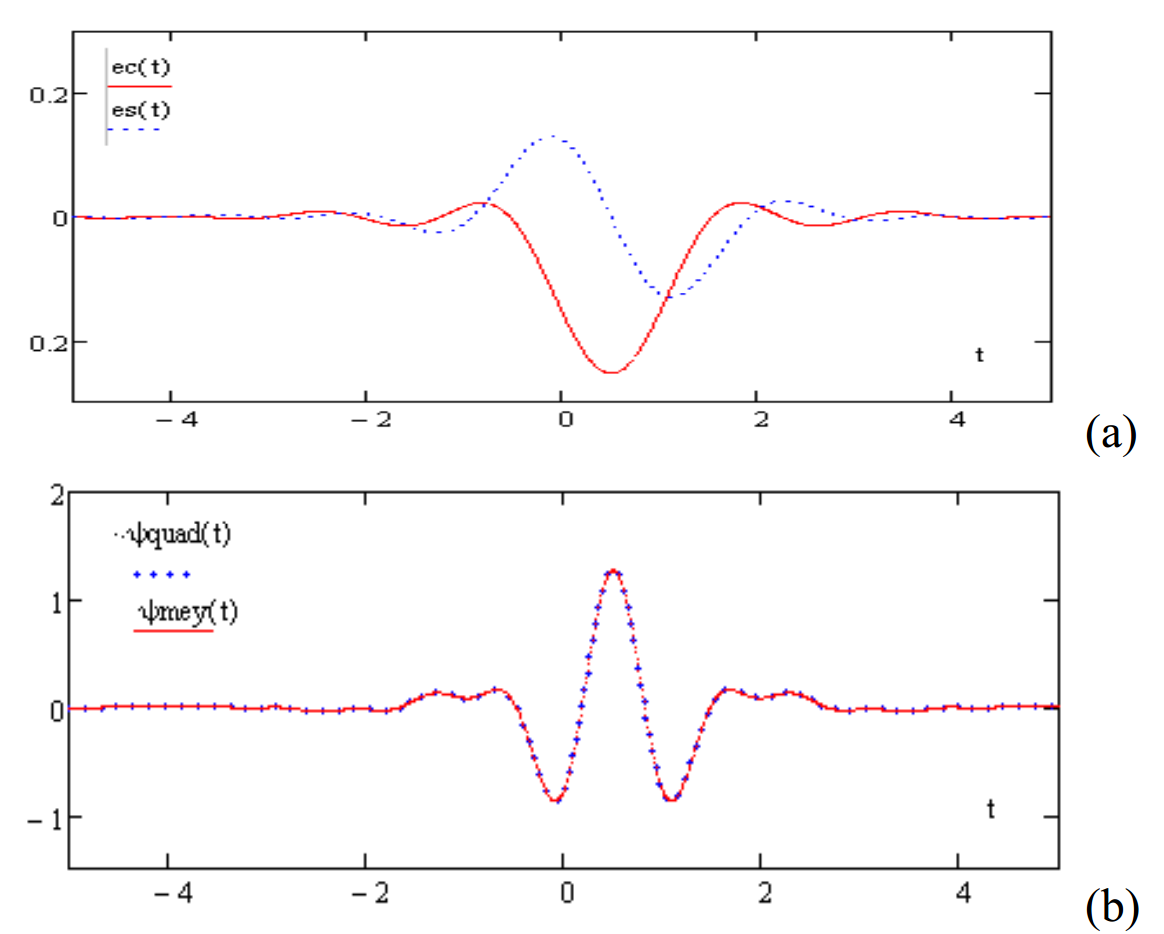}
\captionof{figure}{a) phase and quadrature components of the representation of Meyer wavelet. b) Generating Meyer wavelet from the phase and quadrature components described in (Fig.~\ref{fig:fig5}a), according to the representation bandpass Eqn. (13). The function is obtained by inverse Fourier transform in (11).}
\label{fig:fig5}
\end{figure}

\begin{figure}
\centering
\includegraphics[scale=0.45]{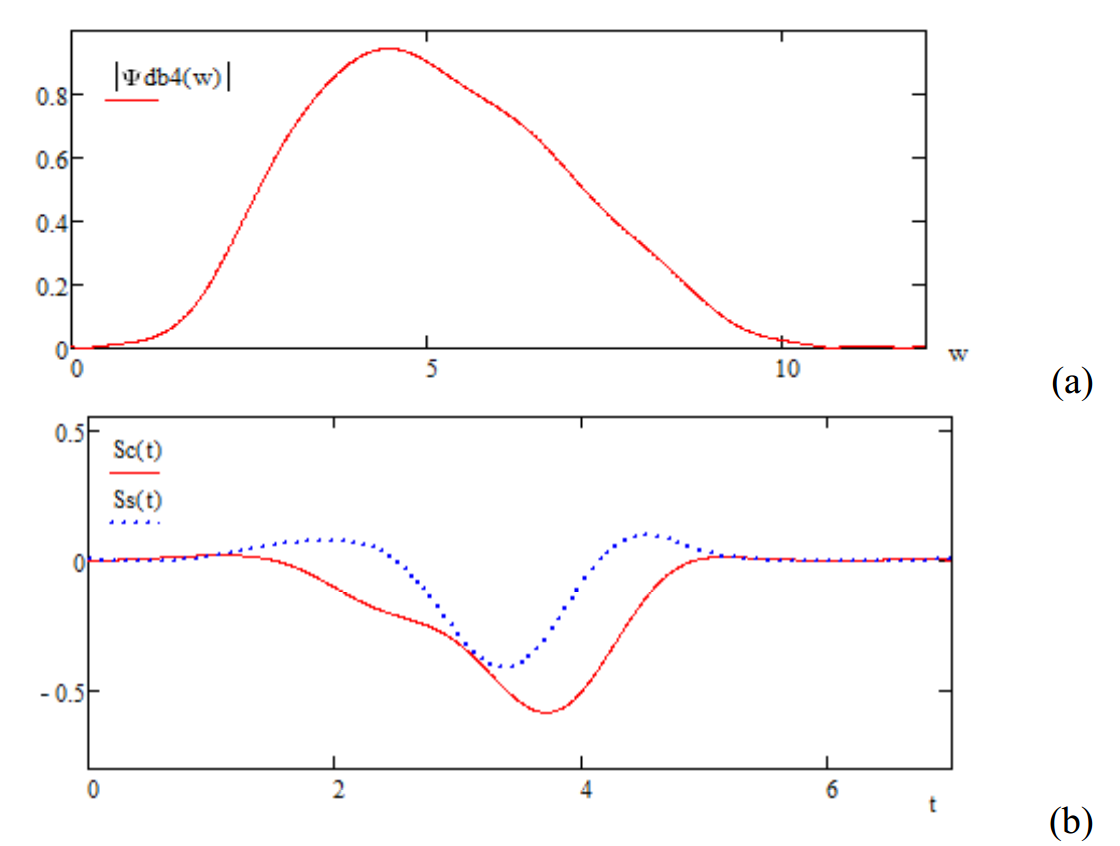}
\captionof{figure}{Kotel'nikov decomposition to db4 wavelet using analytical approaches proposed in \cite{deO4}.}
\label{fig:fig6}
\end{figure}

\begin{figure}
\centering
\includegraphics[scale=0.45]{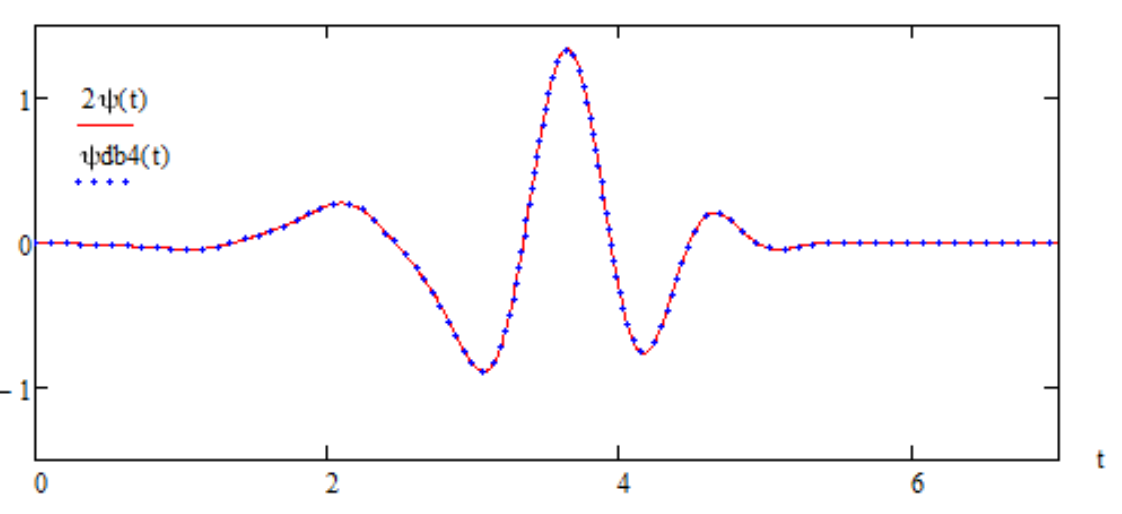}
\captionof{figure}{db4 recovery using in-phase and quadrature components.}
\label{fig:fig7}
\end{figure}

The direct application of Kotel'nikov theorem for a wavelet-mother spectrum (effectively) confined in the spectral range $[f_m, f_M]$ results in the representation:\\
\begin{equation}
\psi(t)=e(t).cos \left [ \pi . (f_M+f_m).t+\theta (t) \right ],
\end{equation}
$e(t)$ and $\theta (t)$ being baseband processes, spectra limited in $[0,B/2]$, and $B:=(f_M-f_m)$ the bandwidth of the wavelet. An alternative way is to consider a complex envelope, modulated by a carrier:\\
\begin{equation}
\psi(t):=\Re e \left \{ S_b(t).e^{j \pi (f_M+f_m) t} \right \},
\end{equation}
where $S_b(t):=e(t).e^{j \theta(t)}$ is a (complex) baseband signal. Directly, as in signal analysis in Telecommunication systems \cite{Vetterli}, \cite{Taub}, it is investigated as a description of the components processes, $e(t)$ and $\theta (t)$, from the waveform $\psi(t)$. Despite the fact that the envelope can be extracted from the use of classical envelope detector \cite{Vetterli}, we opted for an approach to an analytical formulation for both processes, envelopment and phase. Thus, we deal with synchronous detection \cite{Vetterli} to ``demodulate'' the wavelet. The frequency of virtual carrier is exactly the central point of the mother-wavelet spectrum:
\begin{equation}
w_c :=\frac{w_m+w_M}{2} \Rightarrow f_c :=\frac{w_m+w_M}{2} ~~(Hz).
\end{equation}
The detection is made in ``in phase'' and ``in quadrature'' components, using low-pass filters for the respective components. Fig.~\ref{fig:fig1} illustrates the decomposition process. The analysis is conveniently separated into upper and lower branches.\\
\\
Fig.~\ref{fig:fig1} shows two analysis branches:
\begin{subequations}
\begin{align}
a)~upper~branch \nonumber\\
\psi_c(t) & :=\psi(t).cos\left [ \frac{w_m+w_M}{2}t \right ], \\
b)~lower~Branch \nonumber\\
\psi_s(t) & :=\psi(t).sin\left [ \frac{w_m+w_M}{2}t \right ].
\end{align}
\end{subequations}
Using the representation of Kotel'nikov in (1) and substituting in (4a), we obtain the in-phase component:
\begin{equation}
\psi_c(t)=\frac{e(t).cos\left ( \theta (t) \right )}{2}+\frac{e(t)}{2}.cos\left [ \left ( w_m+w_M \right )t+\theta (t) \right ].
\end{equation}
And the output is low-pass filtered with a cutoff frequency $B/2$, resulting therefore in:
\begin{equation}
S_c(t):={\psi (t)} {|}_{LPF_{ed}} = {}\frac{e(t).cos \left [ \theta (t) \right ]}{2}.
\end{equation}
Similarly, the equations come:
\begin{equation}
\psi_s(t)=-\frac{e(t).sin\left [ \theta (t) \right ]}{2}+\frac{e(t)}{2}.sin\left [ \left ( w_m+w_M \right ) t +\theta (t)\right ]
\end{equation}
And the output is low-pass filtered with a cutoff frequency $B/2$, resulting therefore in:
\begin{equation}
S_s(t):=\psi_s(t) {|}_{LPF_{ed}} =-{}\frac{e(t).sin \left [ \theta (t) \right ]}{2}.
\end{equation}
From the Eqns (6) and (8), derive similar relationships as ``classical equations'' telecommunications theory:
\begin{equation}
e(t)=2.\sqrt{S_{c}^{2}(t)+S_{s}^{2}(t)}  e \theta (t)=tg^{-1}\left ( \frac{S_s(t)}{S_c(t)} \right ).
\end{equation}
The envelope of representation bandpass wavelet has a relationship with the associated scaling function. The corresponding analysis filters to the ``ideal bandpass'' defines a wavelet decomposition using wavelet known as the Shannon, whose spectrum is shown in Fig.~\ref{fig:fig2}.
\begin{equation}
\Psi_{SHA}(w)=\Pi \left ( \frac{w-3 \pi /2}{\pi} \right )+\Pi \left ( \frac{w+3 \pi /2}{\pi} \right ),
\end{equation}
where $\Pi \left ( t \right ):=\left\{\begin{matrix}
1 ~~~if~|t|<1/2\\ 
0~~~otherwise.
\end{matrix}\right.$ is the standard gate function.
The bandpass representation resulting in $S_c(t)=sinc(t)$, where $sinc(t):=sin(\pi t)/\pi t$ whereas $S_s(t)=0$. Note that the symmetry of the spectrum implies a zero quadrature component. This leads to an envelope $e(t)=sinc(t)$, which corresponds exactly to scale with a phase function $\theta (t)0=$ \cite{deO}, \cite{Mallat}. Thus, $\psi_{SHA}(t)=sinc(t).cos(3 \pi t/2)$ a known representation. Many real continuous wavelet are already naturally defined under this representation. Even a wavelet infinite support in both domains, time and frequency, including real Morlet wavelet, $\psi_{Morlet}(t):=e^{-t^2}.cos(\pi \sqrt{2/ln2}.t)$ with $f_0:=1/\sqrt{2ln(2)}$ can be approximated by a limitation of effective support in frequency. It is something in the same token as assuming that speech signals are limited band. At this point, it is worth a discussion like that of Slepian's analysis of bandwidth limitation \cite{Slepian}. Since many wavelets with other spectral asymmetry, the center frequency $w_0$ does not correspond to the middle of the strip as in Eqn(1). Whereas, for example, the Meyer wavelet (Fig.~\ref{fig:fig3}), from which orthogonal wavelets is constructed indefinitely derivable infinite medium (the first non-trivial wavelet introduced) have, in frequency domain, wherein (\cite{Mallat},\cite{Meyer}):
\begin{equation}
\begin{split}
&\Psi_{MEY}(w):= \\
& \left\{\begin{matrix}
 {\frac{1}{\sqrt{2 \pi}} sin\left [ \frac{\pi}{2}\nu \left ( \frac{3|w|}{2 \pi}-1 \right )e^{-jw/2} \right ]}&  2 \pi/3 \leq |w| \leq 4 \pi/3 \\ 
 {\frac{1}{\sqrt{2 \pi}} cos\left [ \frac{\pi}{2}\nu \left ( \frac{3|w|}{4 \pi}-1 \right )e^{-j w/2} \right ]} & 4 \pi/3 \leq |w| \leq 8 \pi/3\\ 
0 & otherwise.
\end{matrix}\right.
\end{split}
\label{eq:maineq}
\end{equation}
where
\begin{equation}
\nu (x):= \left\{\begin{matrix}
0 & \forall x\leq 0 \\ 
x & \forall ~0\leq x\leq 1\\ 
1 & \forall x\geq 1.
\end{matrix}\right.
\end{equation}
For implementation in filter bank, see \cite{Dattorro}. Thus, the representation bandpass results in:
\begin{subequations}
\begin{align}
S_c(w) & =\frac{1}{2} \left ( \Psi_{MEY} \left ( w+\frac{6 \pi}{3} \right )+ \Psi_{MEY} \left ( w-\frac{6 \pi}{3} \right )\right ), \\
S_s(w) & =\frac{j}{2} \left ( \Psi_{MEY} \left ( w+\frac{6 \pi}{3} \right )- \Psi_{MEY} \left ( w-\frac{6 \pi}{3} \right )\right ).
\end{align}
\end{subequations}
The signal spectrum at baseband $\Phi (w)$ is sketched in Fig.~\ref{fig:fig4}, i.e., a wrapper $\Phi(w)=\sqrt{(|S_c(w)|^2+|S_s(w)|^2)}$. 
This corresponds exactly to scale function to the wavelet Meyer. Interestingly, this same expression applied to the wavelet Shannon also results in a scaling function.
Comes to an analytical representation to Meyer wavelet, with components outlined in Fig.~\ref{fig:fig5}a:
\begin{equation}
\psi_{quad}(t)=S_c(t).cos\left ( 2 \pi t \right )+S_s(t).sin\left ( 2 \pi t \right ).
\end{equation}
A similar decomposition can be derived for the db4 wavelet using analytical approaches in quasi-harmonic series proposed in \cite{Vermehren}. The frequency ``central'' $w_0$ decomposition was associated with the peak of the spectrum (Fig.~\ref{fig:fig6}a). The respective components ``in phase'' and ``quadrature'' are shown in Fig. ~\ref{fig:fig6}b. The synthesis of wavelet db4 (approximate) from the two components resulted in exactly the original expression (Fig.~\ref{fig:fig7}), as expected. Another compact support wavelet frequency is the complex ``de Oliveira'' wavelet \cite{deO3}.\\
\begin{figure}
\centering
\includegraphics[scale=2.5]{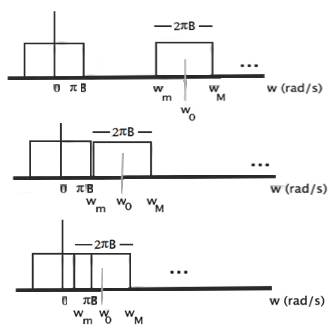}
\captionof{figure}{Analysis of spectral overlap in the decomposition. There are three cases, namely: a) $w_M>3w_M$ no spectral overlap, b) represents the limiting case $w_M=3w_M$, and c) $w_M<3w_M$ occurs spectral overlap in the range $\left ( w_m,\frac{w_M-w_m}{2} \right )$.}
\label{fig:fig8}
\end{figure}

\section{Wavelet Analysis as Filter Bank without Spectral Overlap}

\begin{figure}
\centering
\includegraphics[scale=0.4]{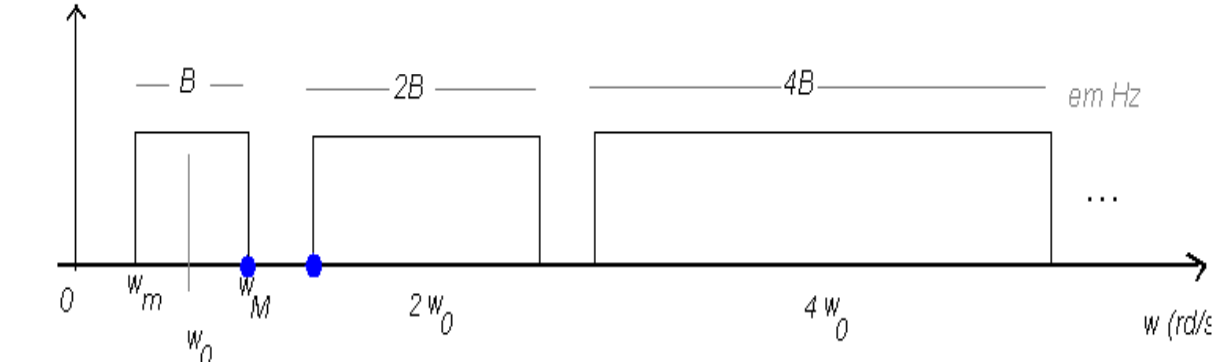}
\captionof{figure}{Dyadic analysis with $Q$-constant filter bank. The band $B$ referenced in the figure is $2 \pi B =w_M - w_m-$(in rd/s).}
\label{fig:fig9}
\end{figure}

\begin{figure}
\centering
\includegraphics[scale=0.44]{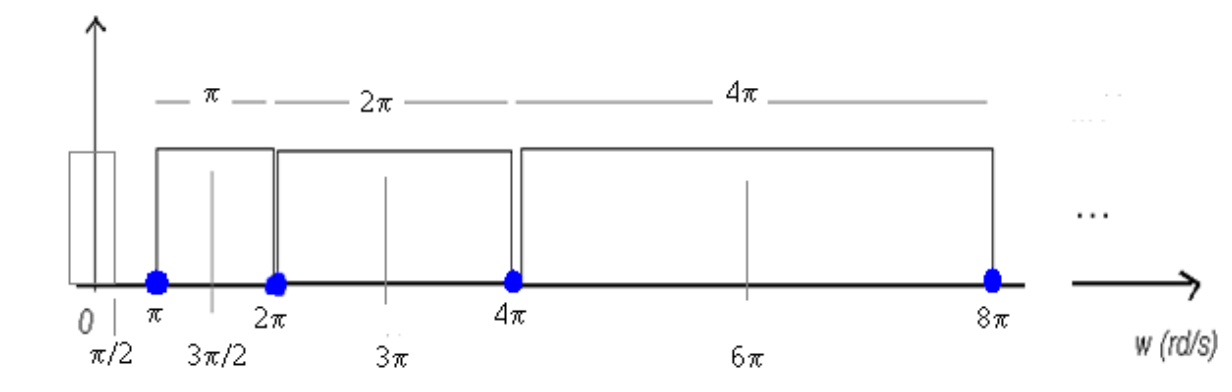}
\captionof{figure}{Orthogonal dyadic filter bank with $Q$-constant in the analysis of the Shannon wavelet. Clearly, there is no spectral overlap.}
\label{fig:fig10}
\end{figure}

\begin{figure}
\centering
\includegraphics[scale=0.44]{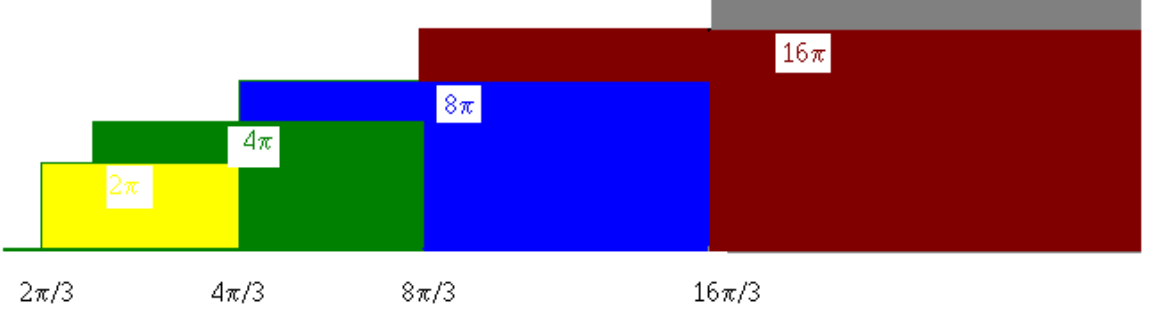}
\captionof{figure}{Filter Bank orthogonal dyadic $Q$-constant in the analysis with Meyer wavelet. Despite the spectral overlap, the analysis is performed correctly.}
\label{fig:fig11}
\end{figure}

The basic properties of wavelet theory are shifting and scaling:
\begin{equation}
\psi_{a,b}(t):=\frac{1}{\sqrt{|a|}}.\psi \left ( \frac{t-b}{a} \right ),~a\neq 0.
\end{equation}
Applying the scaling $a=2$ in the representation shown in (1) to the mother wavelet, is obtained:
\begin{equation}
\begin{split}
e(2t).cos\left [ 2 \pi . (f_M+f_m)t+ \theta(2t) \right ]=\\
\Re e\left \{ S_b (2t).e^{j2 \pi.(f_M+f_m)t}\right \}.
\end{split}
\end{equation}
It is seen that as the frequency is multiplied by $a$, (for $a=2^k$, $k$ integer, the dyadic scale), the complex baseband signal is also scaled in the same ratio, \textit{changing its bandwidth}. This way of interpreting the change of wavelet scale naturally leads to analysis with $Q$-constant. Formatting curve ($shape$) depends on the wavelet, but for overlapping purposes, this is completely irrelevant. Just stick to the sign holder limits, something similar to what occurs in the sampling theorem. By inspection in Fig.~\ref{fig:fig8}, as occurs in the usual demonstration of Theorem of sampling of Nyquist-Shannon-Kotel'nikov \cite{Bissell}, is a condition to ensure that there is no spectral overlap ($aliasing$, the sampling case).\\
\medskip
Enforcing $w_m \geq \pi B$ (or identically $w_0-2 \pi B/2 \geq \pi B$), we arrive at: $f_M \leq 3 f_m$. Otherwise, there is spectral overlap in the range $\left ( f_m,\frac{f_M-f_m}{2} \right )$. For example, to Meyer wavelet, superposition occurs in the interval $\left ( \frac{2 \pi}{3},\frac{8 \pi}{3} \right )$. As for the wavelet Shannon, there is no overlap and the analysis is naturally orthogonal. This result can be summarized in the Proposition 1.
\begin{proposition}
\textit{wavelet parent whose continuous spectrum is essentially confined (effective support) in the band $[f_m,f_M]$ performs an orthogonal analysis provided that $f_M \leq 3f_m.$}\\ 
\end{proposition}
The previous result shows a simple test for orthogonal analysis. In a way, resembles the Nyquist criterion. Results of this type are in the line of thought \textit{Entia non sunt multiplicanda praeter necessitatem} [G. Occam]. To evaluate a orthogonality condition between analyzes at different scales, the filter bank is constructed resulting in Fig.~\ref{fig:fig9}, by successive escalations dyadic. For convenience, we outlined a case without spectral overlap.
Taking some examples of known continuous wavelet, it is seen as a straightforward illustration, the condition of the previous proposition. For the Shannon wavelet, Fig.~\ref{fig:fig10} shows how the analysis is performed. The orthogonality of the filter output follows directly from the Parseval-Plancherel Theorem \cite{deO4}).
It is also noteworthy that the orthogonality established here functions as that obtained in FDM systems \cite{Vetterli}. Best results (more compact wavelets in time) can be found, as occurs in OFDM [16-18] systems, in which, although there is an overlap, the orthogonality condition between the channels remains checked. This is a reason to believe that the condition is found just enough. For Meyer wavelet, the filter bank is shown in Fig.~\ref{fig:fig11}, but note the spectral overlapping in the $[4 \pi /3,8 \pi /3]$ range. The function is chosen such that the combination of superimposed components corresponding to appropriate analysis.
\begin{table}[ht]
\centering
\caption{Frequency bands in Spectral Analysis using Meyer Filter Bank.}
\begin{tabular}{c c c c c }
\hline
\tabularnewline
ref. freq. & spectral range & band & spectral range & bandwidth \tabularnewline
(asymmetric) & (overlapping) & & (FDM bank) & \tabularnewline
\tabularnewline
\hline
\tabularnewline
$4 \pi/3$ & $[2\pi/3,8\pi/3]$ & $2\pi$ & $[4\pi/3,8\pi/3]$& $4\pi/3$\tabularnewline
$8 \pi/3$  & $[4\pi/3,16\pi/3]$ & $4\pi$ & $[8\pi/3,16\pi/3]$& $8\pi/3$\tabularnewline
$16 \pi/3$ & $[8\pi/3,32\pi/3]$ &$8\pi$ &$[16\pi/3,32\pi/3]$ & $16\pi/3$ \tabularnewline
$32 \pi/3$ & $[16\pi/3,64\pi/3]$ &$16\pi$ & $[32\pi/3,64\pi/3]$& $32\pi/3$ \tabularnewline
... &... &... &... &... \tabularnewline
\tabularnewline
\hline
\end{tabular}
\label{table:spectral_range_Meyer}
\end{table}
Note that the overlap occurs between the term $cos(.)$ and $sin(.)$ Stepwise. But the expressions of Meyer wavelet already include different scales for the two cossenoidais terms in Eqn~\ref{eq:maineq}. 
\\
This results in the following expression:
\begin{equation}
\begin{split}
& \Psi(w)= \\
& \left\{\begin{matrix}
\frac{1}{\sqrt{2 \pi}}cas\left [ \frac{\pi}{2}\nu \left ( \frac{3|w|}{4 \pi}-1 \right ) \right ] e^{-jw/2}& 4 \pi/3 \leq |w|\leq 8 \pi /3\\ 
 0 & otherwise,
\end{matrix}\right.
\end{split}
\end{equation}
in which $cas(x):=cos(x)+sin(x)$ is the function cassoidal Hartley. The corresponding spectrum is shown in Fig.~\ref{fig:fig12} below. The corresponding wavelet is shown in Fig.~\ref{fig:fig13}.
\begin{figure}
\centering
\includegraphics[scale=0.49]{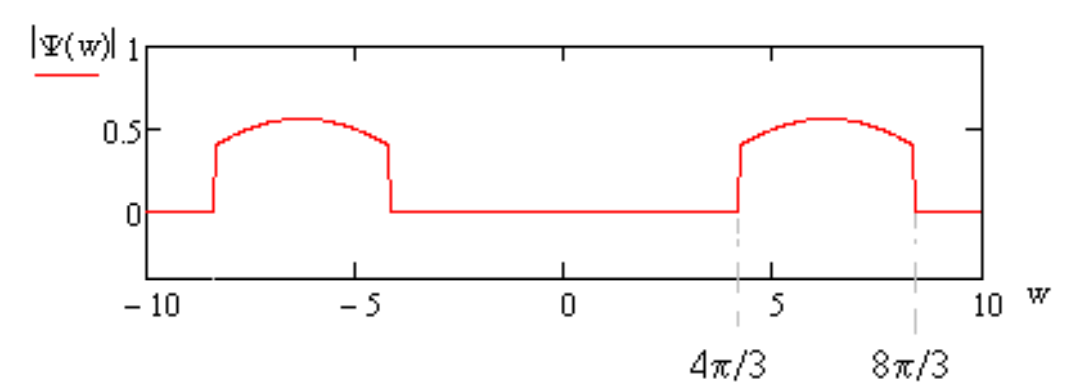}
\captionof{figure}{Analysis for the resulting Meyer bank filters, taking into account the spectral overlap of the side of the bank filters.}
\label{fig:fig12}
\end{figure}
\begin{figure}
\centering
\includegraphics[scale=0.42]{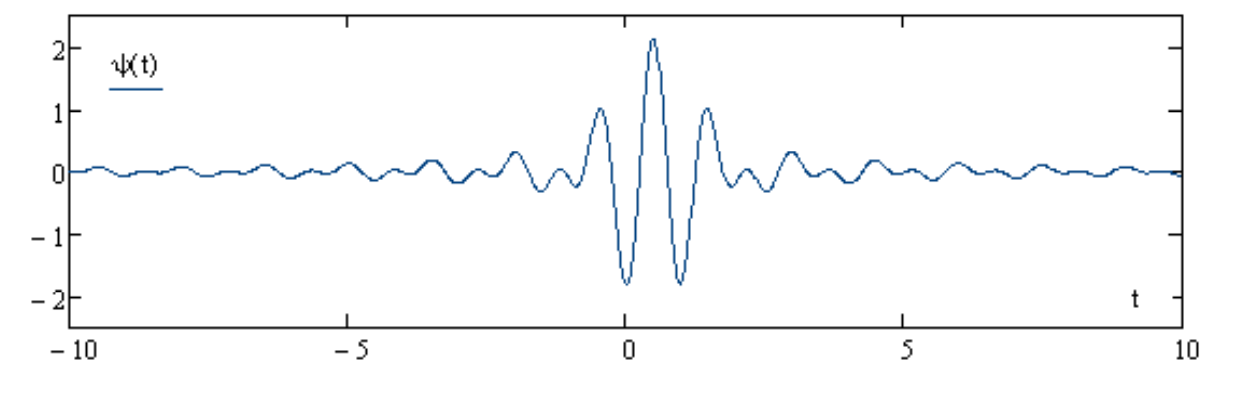}
\captionof{figure}{``equivalent'' Meyer Wavelet, taking into account the effects of superposition of adjacent scales.}
\label{fig:fig13}
\end{figure}
A similar analysis for the ``de Oliveira'' leads to the results shown in Table II. Note that the tracks where there is overlapping spectra between adjacent banks is between $\left [  2 \pi(1-\alpha),2 \pi(1+\alpha)  \right ]$, and the width of the overlap is $4 \pi \alpha$. The corresponding spectrum is sketched in Fig.~\ref{fig:fig14}.\\
\begin{table}[ht]
\centering
\caption{Frequency Bands in Spectral Analysis with de Oliveira filter bank.}
\begin{tabular}{c c}
\hline
\tabularnewline
spectral range & bandwidth\tabularnewline
(overlapping) & \tabularnewline
\tabularnewline
\hline
\tabularnewline
$[\pi(1-\alpha),2\pi(1+\alpha)]$ & $(1+2\alpha) \pi$\tabularnewline
$[2\pi(1-\alpha),4\pi(1+\alpha)]$ & $(1+2\alpha) 2\pi$\tabularnewline
$[4\pi(1-\alpha), 8\pi(1+\alpha)]$& $(1+2\alpha) 4\pi$ \tabularnewline
... &...  \tabularnewline
\tabularnewline
\hline
\tabularnewline
spectral range & bandwidth\tabularnewline
(FDM bank) & \tabularnewline
\tabularnewline
\hline
\tabularnewline
$[ \pi(1+\alpha),2\pi(1+\alpha)]$ & $(1+\alpha)  \pi$\tabularnewline
$[2\pi(1+\alpha),4\pi(1+\alpha)]$ & $(1+\alpha) 2\pi$\tabularnewline
$[4\pi(1+\alpha),8\pi(1+\alpha)]$ & $(1+\alpha) 4\pi$\tabularnewline
... &...  \tabularnewline
\tabularnewline
\hline
\end{tabular}
\label{table:spectral_range_deO}
\end{table}
Note that there is an overlap between the term $cos(.)$ and $cos(.)$-phased. But the expressions of wavelet already include different scales for the two terms in (14), resulting in:
\begin{equation}
\begin{split}
& \Psi(w)= \\
& \left\{\begin{matrix}
0 & |w|\leq \pi(1+\alpha))\\ 
\frac{2}{\sqrt{2 \pi}} & \pi(1+\alpha)\leq |w|\leq 2 \pi(1-\alpha)\\ 
\frac{2}{\sqrt{2 \pi}}.cos\left [ \frac{\left ( |w|-2 \pi(1-\alpha) \right )}{{8 \alpha}}  \right ]  & 2\pi(1-\alpha) \leq |w|\leq 2 \pi(1+\alpha)\\ 
0 & |w|\geq 2 \pi(1+\alpha)
\end{matrix}\right.
\end{split}
\end{equation}
\begin{figure}
\centering
\includegraphics[scale=0.43]{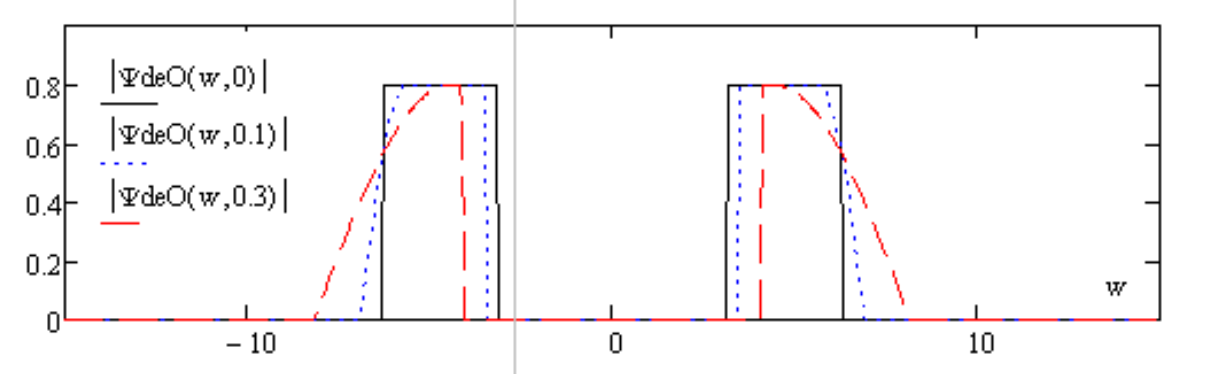}
\captionof{figure}{Analysis for the resulting ``de Oliveira bank filters'', taking into account the spectral overlap of the side of bank filters.}
\label{fig:fig14}
\end{figure}
In Fig.~\ref{fig:fig15}, wavelet waveforms are outlined with different bearing factors corresponding to the analysis filter bank format shown in Fig.~\ref{fig:fig14}.
\begin{figure}
\centering
\includegraphics[height=3cm]{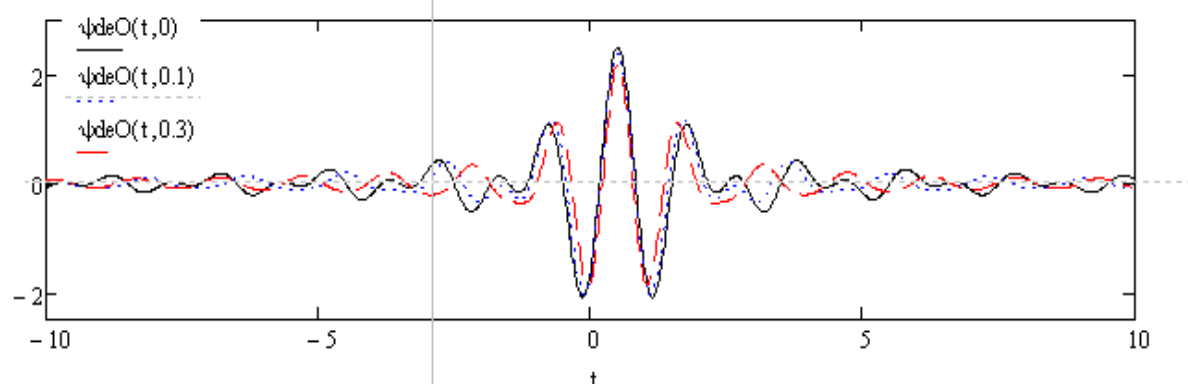}
\captionof{figure}{Real wavelet ``equivalent'' to de Oliveira wavelet, taking into account the effects of superposition of adjacent scales. Three roll-off factor values: $\alpha =0, 0.1$ and $0.3$.}
\label{fig:fig15}
\end{figure}
\section{Conclusions}
Analytical expressions for a bandpass signal representation using a baseband signal are widely spread in telecommunications systems analysis. The procedure for determining the fluctuations of the envelope processes and phase of a wavelet was introduced and illustrated by particular examples of continuous wavelets, including the Daubechies wavelet (db4). Despite the interpretation of wavelets as a filter bank with constant quality factor be old, the approach introduced here can help to better understand the mechanisms involved. \\
The conditions to ensure spectral non-overlap is investigating (the same line as the sampling theorem), which arrives to a sufficient condition on the wavelet mother of spectrum to ensure the spectral orthogonality of the analysis filter bank. Surprisingly, the orthogonal wavelets that result in filter bank where there is some overlap induce an analysis ``equivalent'' with an orthogonal filter bank without spectral overlap. The idea was to combine (principle of superposition) the common areas of the spectrum adjacent filters. As future investigation, it is proposed to try the tools developed in this work in DWT-OFDM wavelet systems \cite{Gupta}.

\section{ACKNOWLEDGMENTS}
The first author thanks D.R. de Oliveira with whom he first shared an early version of this work.

\end{document}